\documentclass{amsart}

\usepackage{hyperref}
\usepackage[T1]{fontenc}

\newtheorem{theo}{Theorem}[section]
\newtheorem{prop}[theo]{Proposition}
\newtheorem{lemm}[theo]{Lemma}

\begin{document}
\title{Quasi-Poisson actions and massive non-rotating BTZ black holes}
\author{S\'ebastien Racani\`ere}
\address{Department of Mathematics, South Kensington Campus\\
Imperial College London\\
SW7 2AZ\\
UK}
\email{s.racaniere@ic.ac.uk}
\urladdr{www.ma.ic.ac.uk/\~{}racani}
\maketitle
\def\g{\mathfrak{g}}
\def\d{\mathfrak{d}}
\def\gsm{{\g_-^\sigma}}
\def\gsp{{\g_+^\sigma}}
\def\Ad{\mbox{\rm Ad}}
\def\id{\mbox{\rm Id}}
\def\im{\mbox{\rm Im}}
\def\SLd{{\mbox{\rm SL}(2,\RR)}}
\def\SLn{{\mbox{\rm SL}(n,\RR)}}
\def\SUd{{\mbox{\rm SU}(2)}}
\def\sin{\mbox{\rm sin}}
\def\cos{\mbox{\rm cos}}
\def\sinh{\mbox{\rm sinh}}
\def\cosh{\mbox{\rm cosh}}
\def\ba{\begin{array}}
\def\ea{\end{array}}
\def\lra{\longrightarrow}
\def\lmt{\longmapsto}
\def\lbr{\lbracket}
\def\rbf{\rbracket}
\let\ldots\dots
\def\F{\mathbf{F}}
\def\FF{\mathbb{F}}
\def\NN{\mathbb{N}}
\def\ZZ{\mathbb{Z}}
\def\QQ{\mathbb{Q}}
\def\RR{\mathbb{R}}
\def\CC{\mathbb{C}}

\begin{abstract}
Using ideas from an article of P.~Bieliavsky, M.~Rooman and Ph.~Spindel on BTZ black holes, I
construct a family of interesting examples of quasi-Poisson actions as defined
by A.~Alekseev and Y.~Kosmann-Schwarzbach. As an application, I obtain a genuine
Poisson structure on $\SLd$ which induces a Poisson structure on a BTZ black hole.
\end{abstract}

\section*{Acknowledgement}
I would like to thank Pierre Bieliavsky for advising me to look at the article~\cite{BRS}
and then answering my questions about it. Thank you to St\'ephane Detournay for talking
to me about physics in a way I could understand and for answering many questions I had about
\cite{BDSR}. And finally, a big thanks to David Iglesias-Ponte for interesting conversations
about \cite{AKS}.

The author was supported by a
Marie Curie Fellowship, EC Contract Ref: HPMF-CT-2002-01850

\section{Introduction}

In \cite{BRS}, P.~Bieliavsky, M.~Rooman and Ph.~Spindel construct a Poisson structure
on massive non-rotating BTZ black holes; in \cite{BDSR}, P.~Bieliavsky, S.~Detournay, Ph.~Spindel
and M.~Rooman construct a star product on the same black hole. The direction of this
deformation is a Poisson bivector field which has the same symplectic leaves as the Poisson
bivector field of \cite{BRS}: roughly speaking, they correspond to orbits under a certain
twisted action by conjugation.

In the present paper, I wish to show how techniques used in \cite{BRS} in conjunction with
techniques of the theory of quasi-Poisson manifolds (see \cite{AKS} and \cite{AKSM}) can be
used to construct an interesting family of manifolds with a quasi-Poisson action and how
a particular case of this family leads to a genuine Poisson structure on a massive non-rotating
BTZ black hole with similar symplectic leaves as in \cite{BRS} and \cite{BDSR}.

\section{Main results}

I will not recall here the basic definitions in the theory of quasi-Poisson manifolds and
quasi-Poisson actions. The reader will find these definitions in A.~Alekseev and
Y.~Kosmann-Schwarzbach~\cite{AKS}, and in A.~Alekseev, Y.~Kosmann-Schwarzbach
and E.~Meinrenken~\cite{AKSM}.

Let $G$ be a connected Lie group of dimension $n$ and $\g$ its Lie algebra, on which $G$ acts
by the adjoint action $\Ad$. Assume we are
given an $\Ad$-invariant non-degenerate bilinear form $K$ on $\g$. For example,
if $G$ is semi-simple, then $K$ could be the Killing form. In the following,
I will denote by $K$ again the linear isomorphism
$$\ba{ccc}
\g & \lra & \g^*\\
x & \lmt & K(x,\cdot).
\ea$$
Let $D=G\times G$ and $\d=\g\oplus\g$ its Lie algebra. Define an $\Ad$-invariant
non-degenerate bilinear form $\langle\,,\rangle$ of signature $(n,n)$ by
$$\ba{ccc}
\d\times\d=(\g\oplus\g)\times(\g\oplus\g) & \lra & \RR\\
((x,y),(x^\prime,y^\prime)) & \lmt & K(x,x^\prime)-K(y,y^\prime).
\ea$$
Assume there is an involution $\sigma$ on $G$ which induces an orthogonal involutive morphism, again
denoted by $\sigma$, on
$\g$. Let $\Delta_+:G\rightarrow D$ and $\Delta_+^\sigma:G\rightarrow D$ be given by
$$\Delta_+(g)=(g,g)$$
and
$$\Delta_+^\sigma(g)=(g,\sigma(g)).$$
Denote by $G_+$ and $G_+^\sigma$ their respective images in $D$. Let $S=D/G_+$ and $S^\sigma=D/G_+^\sigma$. Then both $S$ and $S^\sigma$ are isomorphic
to $G$. The isomorphism between $S$ and $G$ is induced by the map
$$\ba{ccc}
D & \lra & G\\
(g,h) & \lmt & gh^{-1},
\ea$$
whereas the isomorphism between $S^\sigma$ and $G$ is induced by
$$\ba{ccc}
D & \lra & G\\
(g,h) & \lmt & g\sigma(h)^{-1}.
\ea$$
I will use these two isomorphisms to identify $S$ and $G$, and $S^\sigma$ and $G$.
Denote again by $\Delta_+:\g\rightarrow\d$ and $\Delta_+^\sigma:\g\rightarrow\d$ the
morphisms induced by $\Delta_+:G\rightarrow D$ and
$\Delta_+^\sigma:G\rightarrow D$ respectively. Let $\Delta_-:\g\rightarrow\d=\g\oplus\g$ and
$\Delta_-^\sigma:\g\rightarrow\d=\g\oplus\g$ be defined by
$$\Delta_-(x)=(x,-x),$$
and
$$\Delta_-^\sigma(x)=(x,-\sigma(x)).$$
Let $\g_-=\im(\Delta_-)$ and $\g_-^\sigma=\im(\Delta_-^\sigma)$.
We have two quasi-triples $(D,G_+,\g_-)$ and $(D,G_+^\sigma,\g_-^\sigma)$. They induce
two structures of quasi-Poisson Lie group on $D$,
of respective bivector fields $P_D$ and $P_D^\sigma$, and two structures of quasi-Poisson
Lie group on $G_+$ and $G^\sigma_+$ of respective bivector fields $P_{G_+}$ and
$P_{G^\sigma_+}$. I will simply write $G_+$, respectively $G_+^\sigma$, to denote
the group together with its quasi-Poisson structure.
Of course, these quasi-Poisson
structures are pairwise isomorphic. More precisely, the isomorphism
$\id\times\sigma:(g,h)\lmt (g,\sigma(h))$ of $D$ sends $P_D$ on  $P_D^\sigma$
and vice-versa. This isomorphism can be used to deduce some of the results
given at the beginning of the present article from the results of Alekseev and
Kosmann-Schwarzbach~\cite{AKS}; but it takes just as long to redo the computations,
and that is what I do here.

According to \cite{AKS}, the bivector field $P_D$, respectively $P_D^\sigma$,
is projectable onto $S$, respectively $S^\sigma$. Let $P_S$ and $P_{S^\sigma}^\sigma$
be their respective projections. Using the identifications between $S$ and
$G$, and $S^\sigma$ and $G$, one can check that $P_S$ and $P_{S^\sigma}^\sigma$
are the same bivector fields on $G$. What is more interesting, and what I
will prove, is the following Theorem.
\begin{theo}\label{theo:InterestingPoissonAction}
The bivector field $P_D^\sigma$ is projectable onto $S$. Let $P^\sigma_S$ be
its projection. Identify $S$ with $G$ and trivialise their tangent space
using right translations, then for $s$ in $S$ and $\xi$ in $\g^*\simeq T_s^*S$
there is the following explicit formula
$$P_S^\sigma(s)(\xi)=\frac{1}{2}(Ad_{\sigma(s)^{-1}}-Ad_s)\circ\sigma\circ
K^{-1}(\xi).$$
Moreover, the action
$$\ba{ccc}
G^\sigma_+\times S & \lra & S\\
(g,s) & \lmt & gs\sigma(g)^{-1}
\ea$$
of $G_+^\sigma$ on $(S,P_S^\sigma)$ is quasi-Poisson in
the sense of Alekseev and Kosmann-Schwarzbach \cite{AKS}. The image of $P_S^\sigma(s)$,
seen as a map $T^*_s S\lra T_s S$, is tangent to the orbit through $s$ of
the action of $G^\sigma$ on $S$.
\end{theo}

In the setting of the above Theorem, the bivector field $P_S^\sigma$ is $G^\sigma$
invariant; hence if $F$ is a subgroup of $G^\sigma$ and $\bf{I}$ is an $F$-invariant
open subset of $S$ such that the action of $F$ on $\bf{I}$ is principal
then $F\backslash \bf{I}$
is a smooth manifold and $P_S^\sigma$ descends to a bivector field on it.
An application of this remark is the following Theorem.

\begin{theo}\label{theo:TheCaseSL2}
Let $G=\SLd$. Let
$$H=\left[\ba{cc}1&0\\0&-1\ea\right]$$
and choose $\sigma=\Ad_H$. Let
$$\bf{I}=\{\left[\ba{cc}u+x&y+t\\y-t&u-x\ea\right]\mid\,u^2-x^2-y^2+t^2=1, t^2-y^2>0\}$$
be an open subset of $S$. Let $F$ be the following subgroup of $G$
$$F=\{\exp(n\pi H), n\in\NN\}.$$
The quotient $F\backslash \bf{I}$ (together with an appropriate
metric) is a model of massive non-rotating BTZ black
hole (see \cite{BRS}). The bivector field it inherits following the
above remark, is Poisson. Its symplectic leaves consist of the projection
to $F\backslash \bf{I}$ of the orbits of the action of $G^\sigma$ on $S$
except along the
projection of the orbit of the identity. Along this orbit, the bivector field
vanishes and each point forms a symplectic leaf.

In the coordinates $(46)$ of \cite{BRS} (or (\ref{equ:coordinatesonBTZ}) of the present article),
the Poisson bivector field is
\begin{equation}\label{form:PssInCoordinates}
2\cosh^2(\frac{\rho}{2})\sin(\tau)\sinh(\rho) \partial_\tau \wedge \partial_\theta.
\end{equation}
\end{theo}
The above Poisson bivector field should be compared with the one defined
in \cite{BRS} and given by
$$\frac{1}{\cosh^2(\frac{\rho}{2})\sin(\tau)} \partial_\tau \wedge \partial_\theta.$$
The symplectic leaves of this Poisson structure are the images under the projection
$\bf{I}\lra F\backslash\bf{I}$ of the action of $G^\sigma_+$ on $S$.

\section{Let the computations begin}
Throughout the present article, I will use the notations introduced in the previous Section.
To begin with, I will prove that $(D,G^\sigma_+,\g_-^\sigma)$ does indeed form a quasi-triple.

Because $\d=\g\oplus\g$, one also has a decomposition $\d^*=\g^*\oplus\g^*$. One also has
$\d=\gsp\oplus\gsm$ and accordingly $\d^*=\gsp^*\oplus\gsm^*$. Denote $p_\gsp$ and
$p_\gsm$ the projections on respectively $\gsp$ and $\gsm$ induced by the decomposition
$\d=\gsp\oplus\gsm$. So that $1_\d=p_\gsp+p_\gsm$.\\
In this article,
I express results using mostly the
decomposition $\d=\g\oplus\g$. Using it, we have
$$\gsp^*=\{(\xi,\xi\circ\sigma)\mid\,\xi\in\g^*\}$$
and
$$\gsm^*=\{(\xi,-\xi\circ\sigma)\mid\,\xi\in\g^*\}.$$

\begin{prop}
The triple $(D,G^\sigma_+,\g_-^\sigma)$ forms a quasi-triple in the sense of
\cite{AKS}. The characteristic
elements of this quasi-triple as defined in \cite{AKS} and hereby
denoted by $j$, $\F^\sigma$,
$\varphi^\sigma$ and the r-matrix $r_\d^\sigma$ are
$$\ba{ccc}
j:\gsp^* & \lra & \gsm\\
(\xi,\xi\circ\sigma) & \lmt & \Delta_-^\sigma\circ K^{-1}(\xi),
\ea$$
and
$$F^\sigma = 0,$$
and
$$\ba{ccl}
\varphi^\sigma:\bigwedge^3\gsp^* & \lra & \RR\\
((\xi,\sigma\circ\xi),(\eta,\sigma\circ\eta),(\nu,\sigma\circ\eta)) &\lmt& 2K(K^{-1}(\nu),[K^{-1}(\xi),K^{-1}(\eta)]),
\ea$$
and finally the r-matrix
$$\ba{ccc}
r_\d^\sigma:\g^*\oplus\g^* & \lra & \g\oplus\g\\
(\xi,\eta) & \lmt & \frac{1}{2}\Delta_-^\sigma\circ K^{-1}(\xi+\eta\circ\sigma).
\ea$$
\end{prop}
\begin{proof}
It is straightforward to prove that $\d=\g_+^\sigma\oplus\g_-^\sigma$ and that both $\g_+^\sigma$ and
$\g_-^\sigma$ are isotropic in $(\d,\langle\,\rangle)$. This proves that $(D,G^\sigma,\g_-^\sigma)$ is a
quasi-triple.

For $(\xi,\xi\circ\sigma)$ in $\gsp^*$ and $(x,\sigma(x))$ in $\gsp$
$$\langle j(\xi,\xi\circ\sigma),(x,\sigma(x))\rangle=(\xi,\xi\circ\sigma)(x,\sigma(x)).$$
The map $j$ is actually characterised by this last property. The equality
$$\langle \Delta_-^\sigma\circ K^{-1}\circ\Delta_+^{\sigma *}(\xi,\xi\circ\sigma),(x,\sigma(x)\rangle
=(\xi,\xi\circ\sigma)(x,\sigma(x)),$$
proves that
$$\ba{rcl}
j(\xi,\xi\circ\sigma) & = & \Delta_-^\sigma\circ K^{-1}\circ\Delta_+^{\sigma *}(\xi,\xi\circ\sigma)\\
 & = & \Delta_-^\sigma\circ K^{-1}(\xi).
\ea$$

Since $\sigma$ is a Lie algebra morphism, we have
$[\g_-^\sigma,\g_-^\sigma]\subset \g_+^\sigma$. This proves that
$F^\sigma:\bigwedge^2\gsp^*\lra\gsm$, given by
$$F^\sigma(\xi,\eta)=p_\gsm[j(\xi),j(\eta)],$$
vanishes.

I will now compute $\varphi^\sigma$. It is defined as
$$\ba{rcl}
\varphi^\sigma((\xi,\sigma\circ\xi),(\eta,\sigma\circ\eta),(\nu,\sigma\circ\nu)) &=& (\nu,\sigma\circ\nu)\circ
p_\gsp([j(\xi,\sigma\circ\xi),j(\eta,\sigma\circ\eta)])\\
 &=& \langle j(\nu,\sigma\circ\nu),[j(\xi,\sigma\circ\xi),j(\eta,\eta\circ\eta)]\rangle\\
 &=& \langle \Delta_-^\sigma\circ K^{-1}(\nu),[\Delta_-^\sigma\circ K^{-1}(\xi),\Delta_-^\sigma\circ K^{-1}
(\eta)\rangle\\
 &=& 2K(K^{-1}(\nu),[K^{-1}(\xi),K^{-1}(\eta)]).
\ea$$
Finally, the r-matrix is defined as
$$\ba{ccc}
r_\d^\sigma:\gsp^*\oplus\gsm^* & \lra & \gsp\oplus\gsm\\
((\xi,\xi\circ\sigma),(\eta,\eta\circ\sigma)) & \lmt & (0,j(\xi,\xi\circ\sigma)).
\ea$$
If $(\xi,\eta)$ is in $\d^*=\g^*\oplus\g^*$ then its decomposition in $\gsp^*\oplus\gsm^*$ is
$\frac{1}{2}((\xi+\eta\circ\sigma,\xi\circ\sigma+\eta),(\xi-\eta\circ\sigma,-\xi\circ\sigma+\eta))$. The result follows.
\end{proof}

I now wish to compute the bivector $P_D^\sigma$ on $D$. By definition, it is equal to $(r_\d^\sigma)^\lambda-
(r_\d^\sigma)^\rho$, where the upper script $\lambda$ means the left invariant section of $\Gamma(T D\otimes T D)$
generated by $r_\d^\sigma$, while the upper script $\rho$ means the right invariant section of $\Gamma(T D\otimes T D)$
generated by $r_\d^\sigma$.

\begin{prop}
Identify $T_d D$ to $\d$ by right translations. The value of $P_D^\sigma$
at $d=(a,b)$ is
$$\ba{ccc}
\d^*=\g^*\oplus\g^* & \lra & \d=\g\oplus\g\\
(\xi,\eta) & \lmt & \frac{1}{2}(K^{-1}(\eta\circ\sigma\circ(\Ad_{\sigma(b)a^{-1}}-1)),-K^{-1}
(\xi\circ\sigma\circ(\Ad_{\sigma(a)b^{-1}}))).
\ea$$
\end{prop}
\begin{proof}
Fix $d=(a,b)$ in $D$. I choose to trivialise the tangent bundle, and its
dual, of $D$ by using right translations. See $(r_\d^\sigma)^\rho$
as a map from $T^*D$ to $T D$. If $\alpha$ is in $\d^*$, then
$$(r_\d^\sigma)^\rho(d)(\alpha^\rho)=(r_\d^\sigma(\alpha))^\rho(d),$$
whereas
$$(r_\d^\sigma)^\lambda(d)(\alpha^\rho)=(\Ad_d\circ r_\d^\sigma(\alpha\circ\Ad_d))^\rho(d).$$
Thus $P_D^\sigma$ at the point $d=(a,b)$ is
$$\ba{ccc}
\d^*=\g^*\oplus\g^* & \lra & \d=\g\oplus\g\\
(\xi,\eta) & \lmt & -\frac{1}{2}\Delta_-\circ K^{-1}(\xi+\eta\circ\sigma)+\frac{1}{2}\Ad_d\circ
\Delta_-\circ K^{-1}(\xi\circ\Ad_a+\eta\circ\Ad_b\circ\sigma).
\ea$$
The above description of $P_D^\sigma$ can be simplified:
$$\ba{cl}
  & P_D^\sigma(d)(\xi,\eta)\\
= & -\frac{1}{2}\Delta_-\circ K^{-1}(\xi+\eta\circ\sigma)+\\
  & \frac{1}{2}(\Ad_a\circ K^{-1}(\xi\circ\Ad_a+\eta\circ\Ad_b\sigma),
 -\sigma\circ\Ad_{\sigma(b)}\circ K^{-1}(\xi\circ\Ad_a+\eta\circ\Ad_b\circ\sigma))\\
= & -\frac{1}{2}\Delta_-\circ K^{-1}(\xi+\eta\circ\sigma)+\\
  & \frac{1}{2}( K^{-1}(\xi+\eta\circ\Ad_b\sigma\circ\Ad_{a^{-1}}),
 -\sigma\circ K^{-1}(\xi\circ\Ad_{a\sigma(b)^{-1}}+\eta\circ\Ad_b\circ\sigma\circ\Ad_{\sigma(b)^{-1}}))\\
= & -\frac{1}{2}(K^{-1}(\xi+\eta\circ\sigma),-\sigma\circ K^{-1}(\xi+\eta\circ\sigma))+\\
  & \frac{1}{2}(K^{-1}(\xi+\eta\circ\sigma\circ\Ad_{\sigma(b)a^{-1}}),
 -\sigma\circ K^{-1}(\xi\circ\Ad_{a\sigma(b)^{-1}}+\eta\circ\sigma))\\
= & \frac{1}{2}(K^{-1}(\eta\circ\sigma\circ(\Ad_{\sigma(b)a^{-1}}-1)),-\sigma\circ
K^{-1}(\xi\circ(\Ad_{a\sigma(b)^{-1}}-1)))\\
= & \frac{1}{2}(K^{-1}(\eta\circ\sigma\circ(\Ad_{\sigma(b)a^{-1}}-1)),-
K^{-1}(\xi\circ\sigma\circ(\Ad_{\sigma(a)b^ {-1}}-1))).
\ea$$
\end{proof}

It follows from \cite{AKS} that $P_D^\sigma$ is projectable on $S^\sigma=D/G_+^\sigma$.
Actually, the following is also true
\begin{prop}\label{prop:ProjectionOfPdsigmaOnS}
The bivector $P_D^\sigma$ is projectable to a bivector $P_S^\sigma$ on $S=D/G_+$.
Identify $S$ with $G$ through the map
$$\ba{ccc}
D & \lra & G \\
(a,b) & \lmt & ab^{-1}.
\ea$$
Trivialise the tangent space to $G$, and hence to $S$, by right translations. If $s$ is in $S$, then using the above identification, $P_S^\sigma$ at the point $s$ is
\begin{equation}\label{form:ExplicitPss}
P_S^\sigma(s)(\xi)=\frac{1}{2}(\Ad_{\sigma(s)^{-1}}-\Ad_s)\circ\sigma\circ
K^{-1}(\xi).
\end{equation}
\end{prop}
\begin{proof}
Assume $s$ in $S$ is the image of $(a,b)$ in $D$, that is $s=ab^{-1}$. The
tangent map of
$$\ba{ccc}
D & \lra & G \\
(a,b) & \lmt & ab^{-1}
\ea$$
at $(a,b)$ is
$$\ba{ccc}
p:\d & \lra & \g\\
(x,y) & \lmt & x-\Ad_{ab^{-1}}y.
\ea$$
The dual map of $p$ is
$$\ba{ccc}
p^*:\g^* & \lra & \d^*\\
\xi & \lmt & (\xi,-\xi\circ\Ad_{ab^{-1}})
\ea$$
The bivector $P_D^\sigma$ is projectable onto $S$ if and only if for all
$(a,b)$ in $D$ and $\xi$ in $\g^*$, the expression
$$p(P_D^\sigma(a,b)(p^*\xi))$$
depends only on $s=a b^{-1}$. It will then be equal to $P_S^\sigma(s)(\xi)$.
This expression is equal to
$$\ba{cl}
  & p(P_D^\sigma(a,b)(\xi,-\xi\circ\Ad_{ab^{-1}}))\\
= & \frac{1}{2}p((\Ad_{a\sigma(b)^{-1}}-1)\circ\sigma\circ K^{-1}(-\xi\circ\Ad_{ab^{-1}}),
(1-\Ad_{b\sigma(a)^{-1}})\circ\sigma\circ K^{-1}(\xi))\\
= & \frac{1}{2}p((\Ad_{\sigma(ba^{-1})}-\Ad_{a\sigma(a)^{-1}})
\circ\sigma\circ K^{-1}(\xi),
(1-\Ad_{b\sigma(a)^{-1}})\circ\sigma\circ K^{-1}(\xi))\\
= & \frac{1}{2}(\Ad_{\sigma(ba^{-1})}-\Ad_{a\sigma(a)^{-1}}-\Ad_{ab^{-1}}+
\Ad_{a\sigma(a)^{-1}})\circ\sigma\circ K^{-1}(\xi)\\
= & \frac{1}{2}(\Ad_{\sigma(ba^{-1})}-\Ad_{ab^{-1}})\circ\sigma\circ K^{-1}(\xi).
\ea$$
This both proves that $P_D^\sigma$ is projectable on $S$ and gives a formula for the projected bivector.
\end{proof}

To prove that there exists a quasi-Poisson action of $G^\sigma$ on $(S,P_S^\sigma)$,
I must compute $[P_S^\sigma,P_S^\sigma]$, where $[,]$ is the Schouten-Nijenhuis
bracket on multi-vector fields.
\begin{lemm}
For $x$, $y$ and $z$ in $\g$, let $\xi=K(x)$, $\eta=K(y)$ and $\nu=K(z)$. We have
$$\frac{1}{2}[P_S^\sigma(s),P_S^\sigma(s)](\xi,\eta,\nu)
=\frac{1}{4}K(x,[y,\tau_s(z)]+[\tau_s(y),z]-\tau_s([y,z])),$$
where $\tau_s=\Ad_s\circ\sigma-\sigma\circ\Ad_{s^{-1}}$.
\end{lemm}
\begin{proof}
Let $(a,b)$ in $D$ be such that $s=ab^{-1}$. Let $p$ be as in the proof of
Proposition~\ref{prop:ProjectionOfPdsigmaOnS}.
The bivector $P_S^\sigma(s)$ is $p(P_D^\sigma(a,b))$. Hence,
$$[P_S^\sigma(s),P_S^\sigma(s)]=p([P_D^\sigma(a,b),P_D^\sigma(a,b)]).$$
But it is proved in \cite{AKS} that
$$[P_D^\sigma(a,b),P_D^\sigma(a,b)]=(\varphi^\sigma)^\rho(a,b)-(\varphi^\sigma)^\lambda(a,b).$$
Hence
$$\frac{1}{2}[P_S^\sigma(s),P_S^\sigma(s)]=p((\varphi^\sigma)^\rho(a,b))-
p((\varphi^\sigma)^\lambda(a,b)).$$
Now, it is tedious but straightforward and very similar to the above computations to check that
$$p((\varphi^\sigma)^\rho(a,b))(\xi,\eta,\nu)=
\frac{1}{4}K(x,[y,\tau_s(z)]+[\tau_s(y),z]-\tau_s([y,z])),$$
and
$$p((\varphi^\sigma)^\lambda(a,b))(\xi,\eta,\nu)=0.$$
\end{proof}
The group $D$ acts on $S=D/G_+$ by multiplication on the left.
This action restricts to an action of $G^\sigma_+$ on $S$. Identifying
$G$ and $G^\sigma_+$ via $\Delta_+^\sigma$, this action is
$$\ba{ccc}
G\times S & \lra & S\\
(g,s) & \lmt & gs\sigma(g)^{-1}.
\ea$$
The infinitesimal action of $\g$ at the point $s$ in $S$ reads
$$\ba{ccc}
\g & \lra & T_s S\simeq \g\\
x & \lmt & x-\Ad_s\circ\sigma(x),
\ea$$
with dual map
$$\ba{ccc}
T_s^* S\simeq\g^* & \lra & \g^*\\
\xi & \lmt & \xi-\xi\circ\Ad_s\circ\sigma.
\ea$$
Denote by $(\varphi^\sigma)_S$ the induced trivector field on $S$.
If $\xi$, $\eta$ and $\nu$ are in $\g^*$ then
$$(\varphi^\sigma)_S(s)(\xi,\eta,\nu)=\varphi^\sigma(
\xi-\xi\circ\Ad_s\circ\sigma,\eta-\eta\circ\Ad_s\circ\sigma,\nu-\nu\circ\Ad_s\circ\sigma).$$
Computing the right hand side in the above equality is a simple calculation
which proves the following Lemma.
\begin{lemm}
The bivector field $P_S^\sigma$ and the trivector field $(\varphi^\sigma)_S$ satisfy
$$\frac{1}{2}[P_S^\sigma,P_S^\sigma]=(\varphi^\sigma)_S.$$
\end{lemm}

To prove that the action of $G^\sigma_+$ on $(S,P_S^\sigma)$ is indeed quasi-Poisson,
there only remains to prove that $P_S^\sigma$ is $G^\sigma_+$-invariant.

\begin{lemm}
The bivector field $P_S^\sigma$ is $G^\sigma_+$-invariant.
\end{lemm}
\begin{proof}
Fix $g$ in $G\simeq G^\sigma_+$. Denote $\Sigma_g$ the action of $g$ on $S$.
The tangent map of $\Sigma_g$ at $s\in S$ is
$$\ba{ccc}
T_s S\simeq\g & \lra & T_{g s\sigma(g)^{-1}}S\simeq\g\\
x & \lmt & \Ad_g x.
\ea$$
Also, if $\xi$ is in $\g^*$
$$\ba{rl}
  & P_S^\sigma(gs\sigma(g)^{-1})(\xi)\\
= & \frac{1}{2}(\Ad_{g\sigma(s)^{-1}\sigma(g)^{-1}}-
\Ad_{gs\sigma(g)^{-1}})\circ\sigma\circ K^{-1}(\xi)\\
= & \frac{1}{2}\Ad_g\circ(\Ad_{\sigma(s)^{-1}}-
\Ad_s)\circ\Ad_{\sigma(g)^{-1}}\circ\sigma\circ K^{-1}(\xi)\\
= & \Ad_g(P_S^\sigma(s)(\xi\circ\Ad_g))\\
= & (\Sigma_g)_*(P_S^\sigma)(\Sigma_g(s))(\xi).
\ea$$
\end{proof}
\begin{lemm}\label{lemm:ImageOfPss}
Let $s$ be in $S$. The image of $P_S^\sigma(s)$ is
$$\im P_S^\sigma(s)=\{(1-\Ad_s\circ\sigma)\circ(1+\Ad_s\circ\sigma)(y)\mid\, y\in\g\}.$$
In particular, it is included in the tangent space to the orbit through $s$ of the action of $G^\sigma$.
\end{lemm}
\begin{proof}
The image of $P_S^\sigma(s)$ is by Proposition~\ref{prop:ProjectionOfPdsigmaOnS}
$$\im P_S^\sigma(s)  =  \{(\Ad_{\sigma(s)^{-1}}-\Ad_s)\sigma(x)\mid\,
x\in\g\}.$$
The Lemma follows by setting $x=\Ad_s\circ\sigma(y)=\sigma\circ\Ad_{\sigma(s)}(y)$ and
noticing that $(1-(\Ad_s\circ\sigma)^2)=(1-\Ad_s\circ\sigma)\circ(1+\Ad_s\circ\sigma)$.
\end{proof}

This finishes the proof of Theorem~\ref{theo:InterestingPoissonAction}.

Choose $G$ and $\sigma$ as in Theorem~\ref{theo:TheCaseSL2}. The
trivector field $[P_S^\sigma,P_S^\sigma]$ is tangent to the orbit of
the action of $G^\sigma_+$ on $S$. These orbits are of dimension at
most $2$, therefore the trivector field $[P_S^\sigma,P_S^\sigma]$ vanishes
and $P_S^\sigma$ defines a Poisson structure on $\SLd$ which is invariant
under the action
$$\ba{ccc}
\SLd\times\SLd & \lra & \SLd\\
(g,s) & \lmt & g s\sigma(g)^{-1}.
\ea$$
Lemma~\ref{lemm:ImageOfPss} and a simple computation prove
that along the orbit of the identity,
the bivector field $P_S^\sigma$ vanishes; and that elsewhere, its image
coincides with the tangent space to the orbits of the above action.
Recall that in \cite{BRS}, the domain $\bf{I}$ is given by
\begin{eqnarray}
z(\tau,\theta,\rho)=\left[\ba{cc}\sinh(\frac{\rho}{2})+\cosh(\frac{\rho}{2})\cos(\tau) &
\exp(\theta)\cosh(\frac{\rho}{2})\sin(\tau)\\
-\exp(-\theta)\cosh(\frac{\rho}{2})\sin(\tau) &
-\sinh(\frac{\rho}{2})+\cosh(\frac{\rho}{2})\cos(\tau)\ea\right].\label{equ:coordinatesonBTZ}
\end{eqnarray}
This formula also defines coordinates on $I$.
Using Formula~(\ref{form:ExplicitPss}) of
Proposition~\ref{prop:ProjectionOfPdsigmaOnS}
and a computer,
it is easy to check that $P_S^\sigma$ if indeed given by
Formula~(\ref{form:PssInCoordinates}). This ends the proof of
Theorem~\ref{theo:TheCaseSL2}.

\section{A final remark}
One might ask how different is the quasi-Poisson action of $G_+^\sigma$ on
$(S,P_S^\sigma)$ from the usual quasi-Poisson action of $G_+$ on $(S,P_S)$. For
example, if one takes $G=\SUd$, $H=\left[\ba{cc}1 & 0\\0 & -1\ea\right]$ and
$\sigma=\Ad_H$ then the multiplication on the right in $\SUd$ by
$\left[\ba{cc}i & 0\\0 & -i\ea\right]$ defines an isomorphism between the two
quasi-Poisson actions.

Nevertheless, in the example of Theorem~\ref{theo:TheCaseSL2}, the two structures
are indeed different since for example the action of $\SLd$ on itself by conjugation
has two fixed points whereas the action of $\SLd$ on itself used in Theorem~\ref{theo:TheCaseSL2}
does not have any fixed point.

{\providecommand{\bysame}{\leavevmode ---\ }}

\nobreak

\end{document}